\documentclass[12pt]{article}
\usepackage{amsfonts,amscd,a4}
\usepackage{color}
\usepackage{graphicx}
\usepackage{theorem}
\newtheorem{theo}{Theorem}[section]
\newtheorem{prop}[theo]{Proposition}

\newtheorem{coro}[theo]{Corollary}

{\theorembodyfont{\rm}
\newtheorem{example}[theo]{Example}
\newtheorem{remark}[theo]{Remark}}
\newcommand{\bA}{{\bf A}}
\newcommand{\bB}{{\bf B}}

\newcommand{\bG}{{\bf G}}

\newcommand{\cA}{{\mathcal A}}
\newcommand{\cB}{{\mathcal B}}

\newcommand{\cF}{{\mathcal F}}
\newcommand{\cG}{{\mathcal G}}
\newcommand{\cH}{{\mathcal H}}

\newcommand{\cN}{{\mathcal N}}

\newcommand{\cP}{{\mathcal P}}

\newcommand{\cR}{{\mathcal R}}
\newcommand{\cS}{{\mathcal S}}

\newcommand{\sC}{{\mathbb C}}

\newcommand{\sI}{{\mathbb I}}

\newcommand{\sN}{{\mathbb N}}

\newcommand{\sQ}{{\mathbb Q}}
\newcommand{\sR}{{\mathbb R}}

\newcommand{\sT}{{\mathbb T}}

\newcommand{\sZ}{{\mathbb Z}}
\newcommand{\qed}{\rule{1ex}{1ex}}
\begin{document}
\title{Szeg\"o limit theorems for operators with almost
periodic diagonals}
\author{Steffen Roch, Bernd Silbermann}
\date{Dedicated to Vladimir S. Rabinovich on the occasion of his
65th birthday}
\maketitle
\begin{abstract}
The classical Szeg\"o theorems study the asymptotic behaviour of
the determinants of the finite sections $P_n T(a) P_n$ of Toeplitz
operators, i.e., of operators which have constant entries along
each diagonal. We generalize these results to operators which have
almost periodic functions on their diagonals.
\end{abstract}
\section{Introduction} \label{s1}
This paper deals mainly with operators which are constituted by
Laurent or Toeplitz operators and by band-dominated operators. So
we start with introducing some notations and with recalling some
facts about Toeplitz and band-dominated operators and their finite
sections.
\paragraph{Spaces and projections.} Given a non-empty subset $\sI$
of the set $\sZ$ of the integers, let $l^2(\sI)$ stand for the
Hilbert space of all sequences $(x_n)_{n \in \sI}$ of complex
numbers with $\sum_{n \in \sI} |x_n|^2 < \infty$. We identify
$l^2(\sI)$ with a closed subspace of $l^2(\sZ)$ in the natural
way, and we write $P_\sI$ for the orthogonal projection from
$l^2(\sZ)$ onto $l^2(\sI)$.

The set of the non-negative integers will be denoted by $\sZ^+$,
and we write $P$ in place of $P_\sZ^+$ and $Q$ in place of the
complementary projection $I - P$. Thus, $Q = P_{\sZ^-}$ where
$\sZ^-$ refers to the set of all negative integers.

Further, for each positive integer $n$, set
\[
P_n := P_{\{0, \, 1, \, \ldots, \, n-1\}} \quad \mbox{and} \quad
R_n := P_{\{-n, \, -n+1, \, \ldots, \, n-1\}}.
\]
The projections $R_n$ converge strongly to the identity
operator on $l^2(\sZ)$, and the projections $P_n$ converge
strongly to the identity operator on $l^2(\sZ^+)$ when
considered as acting on $l^2(\sZ^+)$ and to the projection $P$
when considered as acting on $l^2(\sZ)$. 

The $C^*$-algebra of all bounded linear operators on a Hilbert
space $H$ will be denoted by $L(H)$.
\paragraph{Functions and operators.} Let $a \in L^\infty(\sT)$,
the $C^*$-algebra of all essentially bounded measurable functions
on the complex unit circle $\sT$, and let
\[
a_j := \frac{1}{2 \pi} \int_0^{2 \pi} a(e^{it}) e^{-ijt} \, dt.
\]
refer to the $j$th Fourier coefficient of $a$. Then the operator
on $l^2(\sZ)$ given by the matrix representation $(a_{i-j})_{i, \,
j \in \sZ}$ with respect to the standard basis of $l^2(\sZ)$
induces a bounded linear operator $L(a)$ on $l^2(\sZ)$, the
so-called {\em Laurent operator with generating function} $a$. The
operator $T(a) := PL(a)P$ acting on $l^2(\sZ^+)$ is called the
{\em Toeplitz operator with generating function} $a$.

Laurent operators are distinguished by their shift invariance. For
$k \in \sZ$, define the {\em shift operator}
\[
U_k : l^2(\sZ) \to l^2(\sZ), \quad (x_n) \mapsto (y_n) \;
\mbox{with} \; y_n = x_{n-k}.
\]
Then $A \in L(l^2(\sZ))$ is a Laurent operator if and only if
$U_{-k} A U_k = A$ for each $k \in \sZ$.

Further, each function $b \in l^\infty(\sZ)$, the $C^*$-algebra of
all bounded sequences on $\sZ$, induces a {\em multiplication
operator}
\[
aI : l^2(\sZ) \to l^2(\sZ), \quad (x_n) \mapsto (a_n x_n).
\]
Let $X$ be a $C^*$-subalgebra of $L^\infty(\sT)$ and $Y$ be a
shift invariant $C^*$-subalgebra of $l^\infty(\sZ)$. The latter
means that $U_{-k} a \in Y$ whenever $a \in Y$ (here we allow the 
operators $U_{-k}$ to act on $L^\infty(\sZ)$ in the obvious way). 
We let $\cA_{X, \, Y}(\sZ)$ stand for the smallest closed
$C^*$-subalgebra of $L(l^2(\sZ))$ which contains all Laurent
operators $L(a)$ with $a \in X$ and all multiplication operators
$bI$ with $b \in Y$. Similarly, we write $\cA_{X, \, Y}(\sZ^+)$
for the smallest closed $C^*$-subalgebra of $L(l^2(\sZ^+))$ which
contains all Toeplitz operators $T(a)$ with $a \in X$ and all
operators $PbP$ with $b \in Y$. So $\cA_{L^\infty (\sT), \,
\sC}(\sZ)$ is the $C^*$-algebra of all Laurent operators, which is
$^*$-isomorphic to the algebra $L^\infty(\sT)$, and
$\cA_{L^\infty(\sT), \, \sC}(\sZ^+)$ is the smallest closed
subalgebra of $L(l^2(\sZ^+))$ which contains all bounded Toeplitz
operators.

Of particular interest are the algebra $X = C(\sT)$ of the
continuous functions on $\sT$ and the algebra $Y = AP(\sZ)$ of the
almost periodic functions. A function $a \in l^\infty(\sZ)$ is
called {\em almost periodic} if the set of all multiplication
operators $U_{-k} a U_k$ with $k \in \sZ$ is relatively compact in
the norm topology of $L(l^2(\sZ))$ or, equivalently, in the norm
topology of $l^\infty(\sZ)$.

The operators in $\cA_{C(\sT), \, l^\infty(\sZ)}(\sZ)$ are usually
referred to as {\em band-dominated operators}, and the operators
in $\cA_{C(\sT), \, Y}(\sZ)$ are called {\em band-dominated operators
with coefficients in} $Y$. To operators in $\cA_{C(\sT), \,
l^\infty(\sZ)}(\sZ^+)$, we also refer as band-dominated operators
over $\sZ^+$. The reason for the notion {\em band-dominated} is
that continuous functions on $\sT$ can be uniformly approximated
by trigonometric polynomials, hence, operators in $\cA_{C(\sT), \,
l^\infty(\sZ)}$ can be approximated by band operators in the norm
of $L(l^2(\sZ))$. We will usually write $\cA_Y(\sZ)$ and
$\cA_Y(\sZ^+)$ in place of $\cA_{C(\sT), \, Y}(\sZ)$ and
$\cA_{C(\sT), \, Y}(\sZ^+)$, respectively, which is consistent
with the notations in \cite{RRS4,Roc6}. It is easy to see that
$PAP \in \cA_Y(\sZ^+)$ whenever $A \in \cA_Y (\sZ)$.
\paragraph{Szeg\"o theorems for Toeplitz matrices.} There are
several ways to express the so-called first Szeg\"o limit theorem,
and there are several kinds of hypotheses under which the theorem
holds. A version which is convenient for us is via stability of
the finite sections method. The $n$th {\em finite section} of the
operator $A$ is the operator $P_nAP_n$. Unless otherwise stated, we
will consider this operator as acting on $\mbox{im} \, P_n$. Thus,
$P_n A P_n$ can be represented by an $n \times n$ matrix. Instead
of $P_n T(a) P_n$ we will also write $T_n(a)$.

The sequence $(P_n A P_n)_{n \in \sN}$ of the finite sections of an
operator $A$ acting on $l^2(\sZ^+)$ is said to be stable if the matrices
$P_n A P_n$ are invertible for sufficiently large $n$ and if the
norms of their inverses are uniformly bounded.
\begin{theo}[First Szeg\"o limit theorem] \label{t1}
Let $a \in L^\infty(\sT)$ and suppose that the finite sections
sequence $(T_n(a))_{n \in \sN}$ is stable. Then $T(a)$ is
invertible and
\begin{equation} \label{e2}
\lim_{n \to \infty} \frac{\det T_n(a)}{\det T_{n-1} (a)} = G[a]
\end{equation}
where
\[
G[a] :=1/(P_1 T(a)^{-1}P_1)
\]
and, of course, $P_1 T(a)^{-1}P_1$ stands for the $00$th entry of
$T(a)^{-1}$.
\end{theo}
If $a \in L^\infty(\sT)$ is real-valued and $T(a)$ is invertible,
then the (compact) essential range of $a$ is contained in the open
interval $(0, \, \infty)$ by the Hartman-Wintner theorem (see 2.36
in \cite{BSi1} or Theorem 1.27 in \cite{BSi2}). Thus, the function
$a$ has a real-valued logarithm $\log a \in L^\infty(\sT)$, and it
is not hard to show that
\begin{equation} \label{e3}
G[a] = \exp \left( \frac{1}{2 \pi} \int_0^{2 \pi} (\log a)(e^{it})
\, dt \right) = \exp (\log a)_0.
\end{equation}
Szeg\"o \cite{Sze1} proved (\ref{e2}) under the assumptions that $a
\in L^1(\sT)$, $a \ge 0$ and $\log a \in L^1(\sT)$. The following
theorems provide statements about the eigenvalue distribution of
Toeplitz matrices. One has to distinguish between real-valued
generating functions $a$, in which case the function $f$ has to be
merely continuous, whereas in case of arbitrary bounded functions
$a$, one needs holomorphy of $f$.

They can be derived from Szeg\"o's first limit theorem (compare
the proofs of Theorems 5.9 and 5.10 in \cite{BSi2}). Although this
derivation is not without effort, they are also referred to as
First Szeg\"o limit theorems. In the present paper we will call
them the distributive versions of Theorem \ref{t1}.

For each $n \times n$-matrix $B$, let $\lambda_i(B)$ with $i = 1,
\, \ldots, \, n$ refer to the eigenvalues of $B$. The order of
enumeration is not of importance.
\begin{theo}[First Szeg\"o limit theorem, distributive version I]
\label{t4} \mbox{} \\
Let $a \in L^\infty(\sT)$ be a real-valued function, and let $g$
be any continuous function on the convex hull of the essential
range of $a$. Then
\begin{equation} \label{e5}
\lim_{n \to \infty} \frac{g(\lambda_1 (T_n(a)) + \dots +
g(\lambda_n (T_n(a))}{n} \, = \, \frac{1}{2 \pi}
\int_0^{2 \pi} g(a(e^{it})) \, dt.
\end{equation}
\end{theo}
\begin{theo}[First Szeg\"o limit theorem, distributive version II]
\label{t6} \mbox{} \\%
Let $a$ be an arbitrary function in $L^\infty(\sT)$, and let $g$
be analytic on an open neighborhood of the convex hull of the
essential range of $a$. Then $(\ref{e5})$ holds again.
\end{theo}
It is one thing to settle the convergence (\ref{e2}) and another
one to describe the precise asymptotic behaviour of the
determinants $\det T_n(a)$. The latter is the contents of
the so-called strong Szeg\"o limit theorem, proved by Szeg\"o
\cite{Sze2} for positive generating functions with H\"older continuous
derivative. In the formulation below, there occurs an algebra, 
$W^{0, \, 0} \cap B_{2, \, 2}^{1/2, \, 1/2}$, of continuous functions 
on $\sT$ which is defined in \cite{BSi1}, 10.21.
\begin{theo}[Strong Szeg\"o limit theorem] \label{t7}
Let $a \in W^{0, \, 0} \cap B_{2, \, 2}^{1/2, \, 1/2}$ have no zeros on 
$\sT$ and winding number 0 with respect to the origin. Then
\begin{equation} \label{e8}
\lim_{n \to \infty} \frac{\det T_n(a)}{G[a]^n} = E[a]
\end{equation}
where
\begin{equation} \label{e9}
E[a] = \exp \sum_{k=1}^\infty k (\log a)_k \, (\log a)_{-k}.
\end{equation}
\end{theo}
We will not go into the long and rich history of the Szeg\"o limit
theorems here and refer to \cite{BSi1,BSi2} and to Chapter 2 of
\cite{Sim1} instead. Let us only mention that E. Basor, G. Baxter,
A. B\"ottcher, A. Devinatz, T. Ehrhardt, I. Gohberg, I. Feldman,
I. I. Hirschman, M. Kac, M. G. Krein and  H. Widom are among the
main contributors and that 
\cite{Bax1,Dev1,Dev2,Ehr2,GoF1,Hir1,Til1,Wid1,ZaT1} mark some 
milestones in this field.
\paragraph{About this paper.} This paper is devoted to
generalizations of the classical Szeg\"o limit theorems to several
classes of operators with variable coefficients (whereas Toeplitz
and Laurent operators are considered as operators with constant
coefficients). Particular attention is paid to operators with
almost periodic coefficients for which we will obtain the most
satisfying generalizations of Theorems \ref{t1} -- \ref{t6}. These
results will be discussed in Sections \ref{s3} and \ref{s4} below. 
On the other hand, we have to report that the precise asymptotic 
behaviour of the determinants of an operator with almost periodic 
diagonals still remains mysterious for us. Thus, the question of a 
possible generalization of the strong Szeg\"o limit theorem is still 
open (although Torsten Ehrhardt's wonderful paper \cite{Ehr2} seems 
to offer a comfortable way to attack this problem).

We prepare our discussion in Section 2 by recalling some facts about 
algebras generated by sequences of finite sections and about
band-dominated operators and their finite sections. The results 
cited in this section can be found in \cite{RRS5,Roc6}. The concluding
fifth section is devoted to some applications of our general
Szeg\"o limit theorems.   

The asymptotic behaviour of the determinants of a sequence of finite 
sections can and should be considered in two different settings: for 
operators acting on the two-sided infinite sequences (with the Laurent
operators as an example) and for operators on one-sided infinite
sequences (for instance, the Toeplitz operators). In order to make our
results comparable with the classical Szeg\"o theorems, we will focus
our attention on operators on $l^2(\sZ^+)$. But many of the presented 
results have their counterparts in the world of operators on two-sided 
sequences. 
\section{Preliminaries} \label{s2}
\subsection{Algebras related with finite sections} \label{ss21}
Let $\cP := (P_n)_{n \in \sN}$ where the projections $P_n$ are
defined as in the introduction. Write $\cF^\cP$ for the set of
all sequences $(A_n)$ of operators $A_n : \mbox{im} \, P_n \to
\mbox{im} \, P_n$ for which the strong limits
\[
\mbox{s-lim} \, A_n P_n \quad \mbox{and} \quad \mbox{s-lim} \,
A_n^* P_n
\]
exist, and $\cG$ for the subset of $\cF^\cP$ consisting of all 
sequences $(G_n)$ with $\|G_n\| \to 0$. Provided with the
operations
\begin{equation} \label{e9b}
(A_n) + (B_n) := (A_n + B_n), \quad \lambda(A_n) := (\lambda A_n),
\quad (A_n) \, (B_n) := (A_nB_n),
\end{equation} 
the involution $(A_n)^* := (A_n^*)$ and with the norm
\[
\|(A_n)\| := \sup_{n \in \sN} \|A_n\|,
\]
the set $\cF^\cP$ becomes a $C^*$-algebra, and $\cG$ is a closed
ideal of $\cF^\cP$. We will often use boldface letters to refer to
elements of $\cF^\cP$. For $\bA := (A_n) \in \cF^\cP$, we denote
the strong limit $\mbox{s-lim} \, A_n P_n$ by $W(\bA)$. Thus, $W$
is a $^*$-homomorphism from $\cF^\cP$ onto $L(l^2(\sZ^+))$.

A sequence $(A_n) \in \cF^\cP$ is called stable if the operators
$A_n : \mbox{im} \, P_n \to \mbox{im} \, P_n$ are invertible for
sufficiently large $n$ and if the norms of their inverses are
uniformly bounded. The following simple result is the basis for
the algebraization of several problems from numerical analysis.
\begin{prop}[Kozak] \label{p9a}
A sequence $\bA \in \cF^\cP$ is stable if and only if the coset
$\bA + \cG$ is invertible in the quotient algebra $\cF^\cP/\cG$.
\end{prop}
The spectrum of the coset $\bA + \cG$ in $\cF^\cP/\cG$ will be
denoted by $\sigma_{\cF^\cP/\cG} \, (\bA + \cG)$ or simply by
$\sigma (\bA + \cG)$. It is also called the {\em stability
spectrum} of the sequence $\bA$ and will occur in the formulation
of several results below. Here we only mention the following fact.
\begin{prop} \label{p9b}
Let $\bA = (A_n) \in \cF^\cP$. Then
\[
\sigma_{L(l^2(\sZ^+))} \, (W(\bA)) \subseteq \sigma_{\cF^\cP/\cG}
\, (\bA + \cG),
\]
and for each open neighborhood $U$ of $\sigma_{\cF^\cP/\cG} \,
(\bA + \cG)$ one has
\[
\sigma_{L(\mbox{\scriptsize im} \, P_n)} \, (A_n) \subseteq U
\]
for all sufficiently large $n$.
\end{prop}
The proof of the first assertion is a consequence of Polski's
theorem (Theorem 1.4 in \cite{HRS2}), and the second one follows
easily from the inclusion
\begin{equation} \label{e9c}
\limsup \sigma(A_n) \subseteq \sigma_{\cF^\cP/\cG} \, (\bA + \cG)
\end{equation}
stated in Theorem 3.19 in \cite{HRS2}, where $\limsup$ is the
set-theoretical limes superior. Indeed, suppose there are an open
neighborhood $U$ of $\sigma_{\cF^\cP/\cG} \, (\bA + \cG)$, a
strongly monotonically increasing sequence $\eta : \sN \to \sN$,
and points $\lambda_n \in \sigma(A_{\eta(n)})$ with $\lambda_n
\not\in U$. Since $(A_n)$ is a bounded sequence, the sequence
$(\lambda_n)$ is bounded, too. Hence, it possesses a partial limit
$\lambda^*$ which belongs to $\limsup \sigma(A_n)$ (by definition)
but not to $U$ (since $U$ is open). This contradicts (\ref{e9c}).
\hfill \qed \\[3mm]
It what follows we will have to consider several subalgebras of
$\cF^\cP$. For $X$ and $Y$ as in the introduction, let $\cS_{X, \,
Y}(\sZ^+)$ stand for the smallest closed $C^*$-subalgebra of
$\cF^\cP$ which contains all sequences $(P_nAP_n)$ of finite
sections of operators $A \in \cA_{X, \, Y}(\sZ^+)$. Further we
will often write $\cS_Y(\sZ^+)$ in place of $\cS_{C(\sT), \,
Y}(\sZ^+)$.
\subsection{Band-dominated operators, their Fredholmness and finite
sections} \label{ss22}
Here is a summary of the results from \cite{RRS1} needed in what
follows. A comprehensive treatment of this topic is in
\cite{RRS4}; see also the references mentioned there.
\paragraph{Fredholmness of band-dominated operators.} An operator
$A$ on a Hilbert space $H$ is called Fredholm if both its kernel
$\mbox{ker} \, A := \{x \in H : Ax = 0\}$ and its cokernel
$\mbox{coker} \, A := H/(AH)$ are finite dimensional linear
spaces. There is a Fredholm criterion for a general band-dominated
operator $A$ which expresses the Fredholm property in terms of the
limit operators of $A$. To state this result, we will need a few
notations.

Let $\cH$ stand for the set of all sequences $h : \sN \to \sZ$
which tend to infinity in the sense that given $C > 0$, there is
an $n_0$ such that $|h(n)| > C$ for all $n \ge n_0$. An operator
$A_h \in L(l^2(\sZ))$ is called the {\em limit operator} of $A \in
L(l^2(\sZ))$ with respect to the sequence $h \in \cH$ if
$U_{-h(n)} A U_{h(n)}$ tends $^*$-strongly to $A_h$ as $n \to
\infty$. (By definition, a sequence $(A_n)$ of operators converges 
$^*$-strongly to $A$ if $A_n \to A$ and $A_n^* \to A^*$ strongly.) 
Notice that every operator can possess at most one limit operator 
with respect to a given sequence $h \in \cH$. The set $\sigma_{\! 
op} (A)$ of all limit operators of a given operator $A$ is the 
{\em operator spectrum} of $A$. \index{operator spectrum} 
\index{$\sigma_{op} (A)$}

We write $L^\$ (l^2(\sZ))$ for the set of all operators $A \in
L(l^2(\sZ))$ which own the following compactness property: Every
sequence $h \in \cH$ possesses a subsequence $g$ for which the
limit operator $A_g$ exists. Thus, operators in $L^\$(l^2(\sZ))$
possess, in a sense, {\em many} limit operators. They are also
called operators with {\em rich} operator spectrum (therefore the
notation).
\begin{prop} \label{p10}
$(a)$ $L^\$ (l^2(\sZ))$ is a $C^*$-subalgebra of $L(l^2(\sZ))$.
\\[1mm]
$(b)$ $\cA_{L^\infty(\sT), \, l^\infty(\sZ)} (\sZ) \subseteq L^\$
(l^2(\sZ))$.
\end{prop}
Assertion $(a)$ is Proposition 1.2.6 (a) in \cite{RRS4}. Since
$L^\$ (l^2(\sZ))$ is a closed algebra, assertion $(b)$ will follow
once it has been shown that all bounded Laurent operators and all
bounded multiplication operators belong to $L^\$ (l^2(\sZ))$. The
first inclusion is evident due to the shift invariance of Laurent
operators, and the second one is Theorem 2.1.16 in \cite{RRS4}.
\hfill \qed \\[3mm]
It is not hard to see that every limit operator of a compact
operator is $0$ and that every limit operator of a Fredholm
operator is invertible. A basic result of \cite{RRS1} (see also
Theorems 2.2.1 and 2.5.7 in \cite{RRS4}) claims that the operator
spectrum of a {\em band-dominated operator} is rich enough in
order to guarantee the reverse implications.
\begin{theo} \label{t11}
Let $A \in L(l^2(\sZ))$ be a band-dominated operator. Then the
operator $A$ is Fredholm if and only if each of its limit
operators is invertible and if the norms of their inverses are
uniformly bounded. If $A$ is a band operator, then $A$ is Fredholm
if and only if each of its limit operators is invertible.
\end{theo}
An analogous result holds for band-dominated operators on $\sZ^+$
in which case one has to take into account all limit operators
with respect to sequences $h$ tending to $+ \infty$. (Simply apply
Theorem \ref{t11} to the operator $PAP+Q$, now acting on all of
$\sZ$.) We let $\sigma_\pm(A)$ collect the set of all limit
operators of $A$ which are taken with respect to a sequence
tending to $\pm \infty$.
\paragraph{Finite sections of band-dominated operators.} One way to 
attack stability problems is based on the following observation.
Associate to the sequence $\bA = (A_n) \in \cF^\cP$ the block diagonal
operator
\begin{equation} \label{e72.1} \label{e12}
\mbox{Op} \, (\bA) := \mbox{diag} \, (A_1, \, A_2, \, A_3, \,
\ldots)
\end{equation}
considered as acting on $l^2(\sZ^+)$. It is easy to check that the
sequence $\bA$ is stable if and only if the associated operator
$\mbox{Op} \, (\bA)$ is Fredholm. In general, this stability
criterion seems to be of less use. But if one starts with the
sequence $\bA = (P_n A P_n)$ of the finite sections method of a
band-dominated operator $A$, then one ends up with a
band-dominated operator $\mbox{Op} \, (\bA)$ on $l^2(\sZ^+)$, and
Theorem \ref{t11} applies to study the Fredholmness of $\mbox{Op}
\, (\bA)$. Basically, one has to compute the limit operators of
$\mbox{Op} \, (\bA)$, which leads to the following result (which
is Theorem 3 in \cite{RRS3}). See also Chapter 6 in \cite{RRS4}
and the detailed account on the finite sections method of
band-dominated operators given in \cite{Roc6}.
\begin{theo} \label{t72.4}
Let $A \in L(l^2(\sZ))$ be a band-dominated operator. Then the
finite sections method $(R_nAR_n)_{n \ge 1}$ is stable if and only
if the operator $A$, all operators
\[
Q A_h Q + P \quad \mbox{with} \quad A_h \in \sigma_+(A)
\]
and all operators
\[
P A_h P + Q \quad \mbox{with} \quad A_h \in \sigma_-(A)
\]
are invertible on $l^2(\sZ)$, and if the norms of their inverses
are uniformly bounded. The condition of the uniform boundedness of
the inverses is redundant if $A$ is a band operator.
\end{theo}
If now $A$ is a band-dominated operator on $l^2(\sZ^+)$, then $PAP+Q$
is a band-dominated operator on $l^2(\sZ)$. Moreover, the finite 
sections sequence $(P_n A P_n)$ is stable if and only if the finite
sections sequence $(R_n (PAP+Q) R_n)$ is stable. Specifying Theorem 
\ref{t72.4} to the case of band operators on $l^2(\sZ^+)$ we get the 
following result, where $J$ refers to the unitary operator
\[
l^2(\sZ) \to l^2(\sZ), \quad (Jx)_m := x_{-m-1},
\]
and where we define $\sigma_+(A)$ as $\sigma_+(PAP+Q)$.
\begin{theo} \label{t72.5}
Let $A \in L(l^2(\sZ^+))$ be a band-dominated operator. Then the
finite sections method $(P_nAP_n)_{n \ge 1}$ is stable if and only
if the operator $A$ and all operators
\[
J Q A_h Q J \quad \mbox{with} \quad A_h \in \sigma_+(A)
\]
are invertible on $l^2(\sZ^+)$ and if the norms of their inverses
are uniformly bounded. The condition of the uniform boundedness of
the inverses is redundant if $A$ is a band operator.
\end{theo}
There are generalizations of Theorems \ref{t72.4} and \ref{t72.5}
which can be verified in the same vein as their predecessors. We
mention the result for the finite sections $(P_n A P_n)$ only.
\begin{theo} \label{t72.6}
Let $A \in L(l^2(\sZ^+))$ be a band-dominated operator, and let
$\eta : \sN \to \sN$ be a strongly monotonically increasing
sequence. Then the sequence $(P_{\eta(n)} A P_{\eta(n)})_{n \ge
1}$ is stable if and only if the operator $A$ and all operators $J
Q A_h Q J$ where $A_h$ is a limit operator of $A$ with respect to
a subsequence $h$ of $\eta$ are invertible on $l^2(\sZ^+)$ and if
the norms of their inverses are uniformly bounded. The condition
of the uniform boundedness of the inverses is redundant if $A$ is
a band operator.
\end{theo}
Thus, instead of taking all limit operators of $A$ with respect to
monotonically increasing sequences $h$, one has to consider only
those with respect to subsequences of $\eta$.
\begin{coro} \label{c16}
Let $A \in L(l^2(\sZ^+))$ be a band-dominated operator, and let $h
: \sN \to \sN$ be a strongly monotonically increasing sequence for
which the limit operator $A_h$ exists. Then the sequence
$(P_{h(n)} A P_{h(n)})_{n \ge 1}$ is stable if and only if the
operators $A$ and $J Q A_h Q J$ are invertible.
\end{coro}
\subsection{Band-dominated operators with almost periodic
coefficients} \label{ss23} 
Here we collect some basic facts from \cite{RRS5} which show
that the conclusion of Corollary \ref{c16} can be essentially
simplified if the sequence $h$ is chosen appropriately. 
These results will only be needed in Subsection \ref{ss51} 
(after Theorem \ref{t39k}) below. 

It is one peculiarity of band-dominated operators $A \in 
\cA_{AP} (\sZ)$ that there is a strongly monotonically increasing 
sequence $h : \sN \to \sN$ such that
\begin{equation} \label{e17}
\|U_{-h(n)} A U_{h(n)} - A\| \to 0 \quad \mbox{as} \; n \to
\infty.
\end{equation}
Thus, $A$ is its own limit operator with respect to $h$, and it is
a limit operator in the sense of {\em norm} convergence. We shall
prove this fact in Section \ref{ss53} in a more general context.
Each sequence $h$ with the properties mentioned above is called a 
{\em distinguished sequence for} $A$. If $h$ is a distinguished 
sequence for $A$, then we call $(P_{h(n)} PAP P_{h(n)})$ the 
associated {\em distinguished finite sections method for} $PAP$ 
and $(R_{h(n)} A R_{h(n)})$ the associated {\em distinguished 
finite sections method for} $A$.
\begin{theo} \label{t18}
Let $A \in \cA_{AP}(\sZ)$ and let $h$ be a distinguished sequence
for $A$. Then the sequence $(P_{h(n)} PAP P_{h(n)})$ is stable if
and only if the operators $PAP$ and $JQAQJ$ are invertible.
\end{theo}
Of course, this follows immediately from Corollary \ref{c16}. But
there is also an elementary proof based on (\ref{e17}) which
mimics the proof of the stability of the finite sections method
for invertible Toeplitz operators with continuous generating
function (see \cite{Boe3}, Theorem 4.102 in \cite{PrS1} and
Section 1.4.2 in \cite{HRS2} for the proof in the Toeplitz setting
and \cite{RRS5} for band-dominated operators with almost periodic
coefficients).

It is not always easy to find a distinguished sequence for a given
operator in $\cA_{AP}(\sZ)$. But sometimes it is, and here are two
examples taken from \cite{RRS5}.
\begin{example}[Multiplication operators] \label{ex18a}
For each real number $\alpha \in [0, \, 1)$, the function
\begin{equation} \label{e19}
a : \sZ \to \sC, \quad n \mapsto e^{2 \pi i \alpha n}
\end{equation}
is almost periodic. Indeed, for every integer $k$, $U_{-k} a U_k$
is the operator of multiplication by the function $a_k$ with
$a_k(n) = a(n+k) = e^{2 \pi i \alpha k} a(n)$, i.e.,
\begin{equation} \label{e20}
U_{-k} a U_k = e^{2 \pi i \alpha k} a.
\end{equation}
Let $(U_{-k(n)} a U_{k(n)})$ by any sequence in $\{ U_{-k} a U_k :
k \in \sZ\}$. Due to the compactness of $\sT$, there are a
subsequence $(e^{2 \pi i \alpha k(n(r))})_{r \ge 1}$ of $(e^{2 \pi
i \alpha k(n)})_{n \ge 1}$ and a real number $\beta$ such that
\[
e^{2 \pi i \alpha k(n(r))} \to e^{2 \pi i \beta} \qquad \mbox{as}
\; r \to \infty.
\]
Thus, the functions $a_{k(n(r))} = e^{2 \pi i \alpha k(n(r))} a$
converge in the norm of $l^\infty(\sZ)$ to $e^{2 \pi i \beta} a$, 
whence the almost periodicity of $a$. Thus, every function as in
(\ref{e19}) belongs to $AP(\sZ)$. Conversely, $AP(\sZ)$ is the 
closure in $l^\infty(\sZ)$ of the span of all functions of the form
(\ref{e19}) with $\alpha \in [0, \, 1)$ (\cite{Cor1}, Theorems 1.9 
-- 1.11 and Theorem 1.27).

For the operator spectrum of the operator $aI$ one finds
\[
\sigma_{\! op, \, s} (aI) = \sigma_{\! op, \, n} (aI) = \left\{
\begin{array}{lll}
\{ e^{2 \pi i l/q} \, a : l = 1, \, 2, \, \ldots, \, q \} &
\mbox{if} & \alpha = 2 p/q \in \sQ, \\[1mm]
\{ e^{it} \, a : t \in \sR \} & \mbox{if} & \alpha \not \in \sQ,
\end{array} \right.
\]
Here, $p$ and $q$ are relatively prime integers with $q > 0$.
Indeed, the inclusion $\subseteq$ follows immediately from
(\ref{e20}). The reverse inclusion is evident in case $\alpha \in
\sQ$. If $\alpha \not\in \sQ$, then it follows from a theorem by
Kronecker which states that the set of all numbers $e^{2 \pi i 
\alpha k}$ with integer $k$ lies dense in the unit circle $\sT$.

Next we are looking for distinguished sequences for the operator
of multiplication by the sequence $\ref{e19}$. From (\ref{e20})
we infer that a sequence $h$ is distinguished for $aI$ if and only 
if
\[
\lim_{n \to \infty} e^{2 \pi i \alpha h(n)} = 1
\]
In case $\alpha = p/q  \in \sQ$, the sequence $a$ is $q$-periodic.
Thus, $h(n) := qn$ is a distinguished sequence for $aI$. For non-rational 
$\alpha \in (0, \, 1)$, expand $\alpha$ into its continued fraction
\[
\alpha = \lim_{n \to \infty} \frac{1}{\displaystyle b_1 +
\frac{1}{\displaystyle b_2 + \frac{1}{\displaystyle
\begin{array}{cc} \displaystyle
\ddots & \\
& \displaystyle b_{n-1} + \frac{1}{b_n}
\end{array}}}}
\]
with uniquely determined positive integers $b_i$. Write this
continued fraction as $p_n/q_n$ with positive and relatively prime
integers $p_n, \, q_n$. These integers satisfy the recursions
\begin{equation} \label{e21}
p_n = a_n p_{n-1} + p_{n-2}, \qquad q_n = a_n q_{n-1} + q_{n-2}
\end{equation}
with $p_0 = 0, \, p_1 = 1, \, q_0 = 1$ and $q_1 = a_1$, and one
has for all $n \ge 1$
\begin{equation} \label{e21a}
\left| \alpha - \frac{p_n}{q_n} \right| < \frac{1}{q_n q_{n+1}} <
\frac{1}{q_n^2}.
\end{equation}
Thus,
\[
|\alpha q_n - p_n| \le q_n \left| \alpha - \frac{p_n}{q_n} \right|
\le \frac{1}{q_n} \to 0,
\]
whence
\[
e^{2 \pi i \alpha q_n} = e^{2 \pi i (\alpha q_n - p_n)} \to 1.
\]
Since moreover $q_1 < q_2 < \ldots$ due to the recursion
(\ref{e21}), this shows that the sequence $h(n) := q_n$ belongs to
$\cH_{A, \, n}$ and that $A_h = A$, i.e. $h$ is a distinguished
sequence for the operator $aI$ with $a$ as in (\ref{e19}). \hfill
\qed
\end{example}
\begin{example}[Almost Mathieu operators] \label{ex21b}
These are the operators $H_{\alpha, \, \lambda, \, \theta}$
on $l^2(\sZ)$ given by
\[
(H_{\alpha, \, \lambda, \, \theta} x)_n := x_{n+1} + x_{n-1} +
\lambda x_n \cos 2 \pi (n \alpha + \theta)
\]
with real parameters $\alpha, \, \lambda$ and $\theta$. Thus,
$H_{\alpha, \, \lambda, \, \theta}$ is a band operator with almost
periodic coefficients, and
\[
H_{\alpha, \, \lambda, \, \theta} = U_{-1} + U_1 + aI \qquad
\mbox{with} \qquad a(n) = \lambda \cos 2 \pi (n \alpha + \theta).
\]
For a treatment of the spectral theory of almost Mathieu operators
see \cite{Boc1} and the recently published papers \cite{AvJ1,Puig}
where the long-standing {\em Ten Martini problem} is solved.

As in Example \ref{ex18a} one gets
\[
U_{-k} H_{\alpha, \, \lambda, \, \theta} U_k = U_{-1} + U_1 + a_k
I
\]
with
\begin{eqnarray} \label{e22}
a_k(n) & = & a(n+k) \; = \; \lambda \cos 2 \pi ((n + k) \alpha +
\theta) \nonumber \\
& = & \lambda ( \cos 2 \pi (n \alpha + \theta) \cos 2 \pi k \alpha
- \sin 2 \pi (n \alpha + \theta) \sin 2 \pi k \alpha).
\end{eqnarray}
We will only consider the non-periodic case, i.e., we let $\alpha
\in (0, \, 1)$ be irrational. As in the previous example, we write
$\alpha$ as a continued fraction with $n$th approximant $p_n/q_n$
such that (\ref{e21a}) holds. Then
\[
\cos 2 \pi \alpha q_n = \cos 2 \pi (\alpha q_n - p_n) = \cos 2 \pi
q_n (\alpha - p_n/q_n) \to \cos 0 = 1
\]
and, similarly, $\sin 2 \pi \alpha q_n \to 0$. Further we infer
from (\ref{e22}) that
\[
|(a_{q_n} - a)(n)| \le |\lambda| \, |1 - \cos 2 \pi \alpha q_n| +
|\lambda| \, |\sin \pi \alpha q_n|.
\]
Hence, $a_{q_n} \to a$ uniformly. Thus, $h(n) := q_n$ defines a
distinguished sequence for the Almost Mathieu operator $H_{\alpha,
\, \lambda, \, \theta}$. Notice that this sequence depends on the
parameter $\alpha$ only. \hfill \qed
\end{example}
Theorem \ref{t18} implies the following.
\begin{coro} \label{c23}
Let $A := H_{\alpha, \, \lambda, \, \theta}$ be an Almost Mathieu
operator and $h$ a distinguished sequence for $A$. Then the
following conditions are equivalent: \\[1mm]
$(a)$ the distinguished finite sections method $(P_{h(n)} PAP
P_{h(n)})$ for $PAP$ is stable; \\[1mm]
$(b)$ the distinguished finite sections method $(R_{h(n)} A
R_{h(n)})$ for $A$ is stable; \\[1mm]
$(c)$ the operators $PAP$ and $QAQ$ are invertible.
\end{coro}
If $\theta = 0$, then the Almost Mathieu operator $A = H_{\alpha,
\, \lambda, \, 0}$ is flip invariant, i.e., $J A J = A$. So we
observe in this case that the third condition in Corollary
\ref{c23} is equivalent to the invertibility of $PAP$ alone.

For a different numerical treatment of Almost Mathieu and other
operators in irrational rotation algebras consult \cite{Bro1}.
\section{The first Szeg\"o limit theorem} \label{s3}
\subsection{Operators with rich spectrum}
Let $A$ be an operator on $l^2(\sN)$ for which the finite sections
sequence $(P_nAP_n)$ is stable. Then the matrices $P_nAP_n$ are
invertible for $n$ large enough, and it makes sense to consider
the sequence
\begin{equation} \label{e24}
n \mapsto \frac{\det (P_n A P_n)}{\det (P_{n-1} A P_{n-1})}.
\end{equation}
In case $A = T(a)$ is an invertible Toeplitz operator with
continuous generating function, the sequence (\ref{e24})
converges, and its limit is equal to
\begin{equation} \label{e25}
G[a] := 1/(P_1 T(a)^{-1} P_1)
\end{equation}
by the first Szeg\"o limit theorem \ref{t1}. For general $A$, one
cannot expect convergence of (\ref{e24}) as already the band
operator
\[
A := \mbox{diag} \, \left( \pmatrix{2 & 1 \cr 1 & 2}, \,
\pmatrix{2 & 1 \cr 1 & 2}, \, \pmatrix{2 & 1 \cr 1 & 2}, \,\ldots
\right)
\]
shows. In this case we denote by $\omega(A)$ \index{$\omega(A)$}
the set of all partial limits of the sequence (\ref{e24}). It
turns out that this set can be described via limit operators in
case $A$ is an operator with rich operator spectrum for which the
finite sections method is stable. We prepare the precise statement
of this result by the following proposition. 
\begin{prop} \label{p25a}
Let $A \in L^\$(l^2(\sZ^+))$ be an operator for which the finite
sections sequence $(P_n A P_n)$ is stable, and let $A_h$ be a limit
operator of $A$ with respect to a sequence $h$ tending to $+ \infty$. 
Then the operator $JQA_hQJ$ is invertible on $l^2(\sZ^+)$.
\end{prop}
For band-dominated operators $A$, this result has been already 
stated in Theorem \ref{t72.5}. In fact, it is the elementary of the 
two implications of the equivalence stated in that theorem. It is 
easy to see that this implication holds for arbitrary operators with 
rich spectrum (see also Proposition 1.2.10 in \cite{RRS4}). \hfill 
\qed \\[3mm] 
The previous observation justifies to set (in analogy to (\ref{e25})) 
\begin{equation} \label{e26}
G[A_h] := 1/(P_1 (JQA_hQJ)^{-1} P_1)
\end{equation}
which has to be read as follows: $P_1 (JQA_hQJ)^{-1} P_1$ can be
understood as an $1 \times 1$-matrix, and we identify this matrix
with its only entry, which is a complex number. The fact that this 
number cannot be zero is part of the assertion of the following 
theorem.
\begin{theo} \label{t27}
Let $A \in L^\$(l^2(\sZ^+))$ be an operator for which the finite
sections sequence $(P_n A P_n)$ is stable. Then $P_1 (JQA_hQJ)^{-1} 
P_1 \neq 0$ for all limit operators $A_h$ of $A$, and 
\begin{equation} \label{e28}
\omega(A) = \{ G[A_h] : A_h \in \sigma_+(A) \}
\end{equation}
with $G[A_h]$ defined by $(\ref{e26})$.
\end{theo}
{\bf Proof.} First we show that $P_1 (JQA_hQJ)^{-1} P_1 \neq 0$ for 
every limit operator $A_h$ of $A$. Let $A_h$ be a limit operator of $A$. 
Equivalently, we have to show that the operator
\[
B_1 := P_1 (JQA_hQJ)^{-1} P_1 : \mbox{im} \, P_1 \to \mbox{im} \, P_1 
\]
is invertible. By Kozak's identity (Proposition 7.15 in \cite{BSi1})
this happens if and only if the operator
\[
B_2 := (P - P_1) JQA_hQJ (P - P_1) : \mbox{im} \, (P - P_1) \to \mbox{im} \, 
(P - P_1)
\]
invertible. We multiply the operator $B_2$ from both sides by the flip 
operator $J$ and take into account that $J (P - P_1) J = (I - R_1) Q$ 
to obtain that $B_2$ is invertible if and only if 
\[
B_3 := (I - R_1) Q A_h Q (I - R_1) : \mbox{im} \, Q (I - R_1) \to
\mbox{im} \, Q (I - R_1)
\]
is invertible. Since $U_1 (I - R_1) Q U_{-1} = Q$, the invertibility 
of $B_3$ is equivalent to the invertibility of the shifted operator 
\begin{eqnarray*}
B_4 & := & U_1 B_3 U_{-1} \\
& = & U_1 (I - R_1) Q U_{-1} U_1 A_h U_{-1} U_1 Q (I - R_1) U_{-1} \\
& = & Q U_1 A_h U_{-1} Q : \mbox{im} \, Q \to \mbox{im} \, Q.
\end{eqnarray*}
It is finally obvious that $B_4$ is invertible if and only if the operator
\[
B_5 := Q U_1 A_h U_{-1} Q + P : l^2(\sZ) \to l^2(\sZ)
\]
is invertible. Since $U_1 A_h U_{-1}$ is also a limit operator of $A$ 
(with respect to the sequence $h^\prime (n) := h(n) -1$), the invertibility 
of $B_5$ follows from the stability of the finite section method 
$(P_nAP_n)$ and from Proposition \ref{p25a}. This settles the first assertion 
of the theorem.

For the second assertion, let $n$ a positive integer and consider the 
operators
\begin{equation} \label{e28a}
W_n : l^2(\sZ^+) \to l^2(\sZ^+), \quad (x_0, \, x_1, \, \ldots)
\mapsto (x_{n-1}, \, x_{n-2}, \, \ldots, \, x_0, \, 0, \, 0, \,
\ldots).
\end{equation}
If the finite sections method $(P_nAP_n)$ is stable, then the
operators $W_nAW_n$, considered as acting on $\mbox{im} \, W_n =
\mbox{im} \, P_n$, are invertible for large $n$, and
\[
\frac{\det (P_{n-1}AP_{n-1})}{\det (P_nAP_n)} = \frac{\det
(W_{n-1}AW_{n-1})}{\det (W_nAW_n)} =: \beta_n.
\]
By Cramer's rule, $\beta_n$ equals the first component of the
solution $x^{(n)}$ to the equation
\[
W_n A W_n x^{(n)} = (1, \, 0, \, 0, \, \dots, \, 0)^T.
\]
Let now $\alpha \in \omega(A)$, and let $h : \sN \to \sN$ be a
sequence tending to infinity such that $\alpha^{-1} = \lim
\beta_{h(n)}$. Since $A$ has a rich operator spectrum, there is a
subsequence $g$ of $h$ such that the limit operator
\[
A_g = \mbox{s-lim} \, U_{-g(n)} A U_{g(n)} \in L(l^2(\sZ))
\]
exists. Then also the strong limit on $l^2(\sN)$
\begin{eqnarray*}
\lefteqn{\mbox{s-lim} \, J U_{-g(n)} P_{g(n)} A P_{g(n)} U_{g(n)} J} \\ 
&& = \mbox{s-lim} \, J (U_{-g(n)} P_{g(n)} U_{g(n)}) \,
(U_{-g(n)} A U_{g(n)}) \, (U_{-g(n)} P_{g(n)} U_{g(n)}) J
\end{eqnarray*}
exists and is equal to $JQA_gQJ$. Since $J U_{-n} P_n = W_n$ and
$P_n U_n J = W_n$, this shows that the strong limit $\mbox{s-lim}
\, W_{g(n)} A W_{g(n)}$ exists and that this limit is equal to
$JQA_gQJ \in L(l^2(\sN))$. So one can consider $(W_{g(n)} A
W_{g(n)})_{n \in \sN}$ as a stable and convergent approximation
sequence for the operator $JQA_gQJ$. In particular, the solutions
$x^{(n)}$ to the equation
\begin{equation} \label{e29}
W_{g(n)} A W_{g(n)} x^{(n)} = (1, \, 0, \, 0, \, \dots, \, 0)^T
\end{equation}
converge in the norm of $l^2(\sN)$ to the solution $x$ to the
equation
\begin{equation} \label{e30}
JQA_gQJ x = (1, \, 0, \, 0, \, \dots )^T.
\end{equation}
Thus, the first component $\beta_{g(n)}$ of the solution $x^{(n)}$
to equation (\ref{e29}) converges to the first component of the
solution $x$ to equation (\ref{e30}). Since the latter one is
equal to
\[
P_1 x = P_1 (JQA_gQJ)^{-1} P_1,
\]
we arrive at $\alpha = (P_1 (JQA_gQJ)^{-1} P_1)^{-1} = G[A_g]$.
This settles the inclusion $\subseteq$ in (\ref{e28}). The
reverse inclusion can be proved by similar arguments. \hfill \qed
\subsection{Operators in the Toeplitz algebra}
By Proposition \ref{p10} $(b)$, the assertion of Theorem \ref{t27}
holds in particular for operators in the algebra
$\cA_{L^\infty(\sT), \, l^\infty(\sZ)} (\sZ^+)$ and, thus, for all
band-dominated operators $A \in \cA_{l^\infty(\sZ)} (\sZ^+)$ and
for all operators $A$ in the Toeplitz algebra $\cA_{L^\infty(\sT),
\, \sC} (\sZ^+)$. The statement for band-operators has been
already proved in \cite{Roc6}, Theorem 7.23, whereas the Toeplitz
case was the subject of Section 7.2.3 in \cite{HRS2}. In the
Toeplitz case, one can complete the assertion of Theorem \ref{t27}
essentially. The point is the following observation.
\begin{prop} \label{p31}
Let $A \in \cA_{L^\infty(\sT), \, \sC} (\sZ^+)$. \\[1mm]
$(a)$ Consider $A$ as an operator on $l^2(\sZ)$ which acts as the
zero operator on $l^2$ over the negative integers. Then the
sequence $(U_{-n} A U_n)_{n \in \sN}$ converges $^*$-strongly on
$l^2(\sZ)$. Its limit is a bounded Laurent operator, i.e., it is
of the form $L(a)$ with $a \in L^\infty(\sT)$. \\[1mm]
$(b)$ The sequence $(W_n A W_n)_{n \in \sN}$ converges
$^*$-strongly on $l^2(\sZ^+)$. Its limit is a bounded Toeplitz
operator, i.e., it is of the form $T(b)$ with $b \in L^\infty(\sT)$.
\\[1mm]
Moreover, $b(t) = \tilde{a}(t) := a(1/t)$ a.e. on $\sT$.
\end{prop}
The function $a$ is also called the {\em symbol} of the operator
$A \in \cA_{L^\infty(\sT), \, \sC} (\sZ^+)$. We denote it by
$s_A$.

For a proof of assertion $(a)$, write $T(a)$ as $P L(a) P$.
Clearly, $U_{-n} L(a) U_n = L(a)$, and one easily checks that
$U_{-n} P U_n \to I$ strongly. Thus,
\[
U_{-n} T(a) U_n \to L(a) \quad \mbox{as} \; n \to \infty.
\]
Assertion $(b)$ follows from $(a)$ since
\[
W_n A W_n = JQ U_{-n} A U_n QJ.
\]
For another proof of $(b)$ (and some facts around it) see Sections
4.3.3 and 7.2.3 in \cite{HRS2}. \hfill \qed \\[3mm]
It follows from Proposition \ref{p31} that the only limit operator
at $+\infty$ of $A \in \cA_{L^\infty(\sT), \, \sC} (\sZ^+)$ is the
Laurent operator $L(s_A)$. Hence, the set $\omega(T(a))$ is the
singleton $\{ G[T(\widetilde{s_A})] \}$ in this case, whence the
convergence of the sequence (\ref{e24}) to this value.
\begin{coro} \label{c32}
Let $A \in \cA_{L^\infty(\sT), \, \sC} (\sZ^+)$ be an operator for
which the finite sections sequence $(P_n A P_n)$ is stable. Then
the sequence $(\ref{e24})$ converges, and its limit is
\[
G[T(\widetilde{s_A})] = 1/(P_1 (T(\widetilde{s_A}))^{-1} P_1).
\]
\end{coro}
\begin{coro} \label{c33}
Let $a \in L^\infty (\sT)$ be such that the finite sections
sequence $(P_n A P_n)$ for the Toeplitz operator $A = T(a)$ is
stable. Then the sequence $(\ref{e24})$ converges, and its limit
is
\[
G[T(\tilde{a})] = 1/(P_1 (T(\tilde{a}))^{-1} P_1).
\]
\end{coro}
In order to show that this corollary indeed reproduces the first
Szeg\"o limit theorem \ref{t1} we have to verify that
\begin{equation} \label{e34}
P_1 T(a)^{-1} P_1 = P_1 (T(\tilde{a}))^{-1} P_1.
\end{equation}
Let $C : l^2(\sZ^+) \to l^2(\sZ^+)$ denote the operator of
conjugation $(x_n) \mapsto (\overline{x_n})$ (which is linear over
the field of the real numbers only). One easily checks that
\[
T(\tilde{a}) = C T(a)^* C \quad \mbox{for each function} \; a \in
L^\infty(\sT).
\]
Hence, $T(a)$ is invertible if and only if $T(\tilde{a})$ is
invertible, and if $B$ is the inverse of $T(a)$, then $CB^*C$ is
the inverse of $T(\tilde{a})$. The 00th entries of $B$ and $CB^*C$
coincide obviously, whence (\ref{e34}). \hfill \qed \\[3mm]
There are two obstacles for the application of Corollary
\ref{c33}. The first one concerns the stability of the finite
sections sequence $(P_n T(a) P_n)$ for which there is no general
criterion known. But there are at least special classes of
generating functions $a \in L^\infty(\sT)$ (e.g., piecewise
continuous or piecewise quasicontinuous functions) for which one
knows that the finite sections sequence for the Toeplitz operator
$T(a)$ is stable if and only the operator $T(a)$ is invertible,
and for which effective criteria for the invertibility of $T(a)$
are available. Details can be found in Section IV.3 in
\cite{GoF1}, Section 4.2 in \cite{HRS2} and Section 2.4 in
\cite{BSi2} for Toeplitz operators with piecewise continuous
generating functions and in Chapter 7 in \cite{BSi1} where a heavy
machinery is developed to attack stability problems.

The second point concerns the constant $G[a] = (P_1 T(a)^{-1}
P_1)^{-1}$ for which one wants to have an effective way of
computation. Under suitable assumptions for the generating
function $a$ (e.g., belonging to the Wiener algebra or being
locally sectorial) one can identify the number $G[a]$ with $1/
\exp (\log a)_0$ with $b_0$ referring to the 0th Fourier
coefficient of the function $b$ (details can be found in Section
5.4 of \cite{BSi2}, for example).

The latter observation offers a also way to determine the constant
$G[A_h]$ in some further instances. Recall that a function $b \in 
l^\infty(\sZ)$ is called {\em slowly oscillating} if the difference 
$b(n+1) - b(n)$ tends to zero as $n \to \pm \infty$. Let $A \in 
L(l^2(\sZ^+))$ be a band-dominated operator with slowly oscillating 
coefficients. It is shown in \cite{LRR1} (see also Theorem 2.9 in 
\cite{Roc6}) that the finite sections method for $A$ is stable if 
and only if the operator $A$ is invertible. Moreover, being 
band-dominated, the operator $A$ has a rich operator spectrum by 
Proposition \ref{p10} $(b)$. Thus, every invertible band-dominated 
operator $A$ with slowly oscillating coefficients satisfies the 
assumptions of Theorem \ref{t27}.

Moreover, in the case at hand, all limit operators of $A$ are shift
invariant (Proposition 2.4.1 in \cite{RRS4}); hence, all partial
limits in $\omega(A)$ are of the form $P_1 T(\tilde{a_h})^{-1}
P_1$ with a certain continuous function $a_h$. If, moreover, $A =
\sum a_k V_k$ satisfies the Wiener condition $\sum \|a_k\|_\infty
< \infty$, then all functions $a_h$ belong to the Wiener algebra,
and one has
\[
P_1 T(\tilde{a_h})^{-1} P_1 = P_1 T(a_h)^{-1} P_1 = 1 / \exp (\log
a_h)_0.
\]
\section{Distributive versions of the first Szeg\"o limit theorem}
\label{s4}
The goal of this section is to prove versions of Theorems \ref{t4}
and \ref{t6} for operators in $\cA_{L^\infty(\sT), \,
AP(\sZ)}(\sZ^+)$. For their formulation, we need some
preparations.

It will be convenient to put the proof into some algebraic
framework which has been developed by Arveson, B\'{e}dos, and
SeLegue \cite{Arv3,Arv4,Bed1,Bed2,SeL1} (see also Section 7.2.1 in
\cite{HRS2}) and which we are going to recall first. For the
reader's convenience, we include the proofs.
\subsection{The F{\o}lner algebra} \label{ss41}
For each operator $A \in L(l^2(\sZ^+))$, let $|A|$ denote its
absolute value, i.e., the non-negative square root of $A^*A$.
Let further $\mbox{tr}$ refer to the canonical trace on the 
finite rank/trace class operators on $l^2(\sZ^+)$, and abbreviate 
the sequence $(P_n)$ to $\cP$. Evidently, $\mbox{tr} \, P_n = n$.
\begin{prop} \label{p41}
The set $\mathfrak{F}(\cP)$ of all operators $A \in L(l^2(\sZ^+))$
with
\begin{equation}\label{e42}
\lim_{n \to \infty} \frac{\mbox{\rm tr} \, (|P_nA -
AP_n|)}{\mbox{\rm tr} \, P_n} = 0
\end{equation}
is a $C^*$-subalgebra of $L(l^2(\sZ^+))$.
\end{prop}
We refer to $\mathfrak{F}(\cP)$ as the {\em F{\o}lner algebra}
\index{F{\o}lner algebra} associated with $\cP$. \\[3mm]
{\bf Proof.} Recall that the set $\cN_1 := \{A \in L(l^2(\sZ^+)) :
\mbox{tr} \, (|A|) < \infty \}$ of the trace class operators is a
two-sided (non-closed) ideal of $L(l^2(\sZ^+))$, that the mapping
$A \mapsto \mbox{tr} \, (|A|)$ defines a norm on $\cN_1$ which
makes this set to a Banach space, and that
\begin{equation} \label{e42a}
|\mbox{tr} \, (A)| \le \mbox{tr} \, (|A|),
\end{equation}
\begin{equation} \label{e43}
\mbox{tr} \, (|A + B|) \le \mbox{tr} \, (|A|) + \mbox{tr} \,
(|B|),
\end{equation}
\begin{equation} \label{e44}
\max \, \{ \mbox{tr} \, (|AC|), \, \mbox{tr} \, (|CA|) \} \le
\|C\| \, \mbox{tr} \, (|A|),
\end{equation}
\begin{equation} \label{e45}
\mbox{tr} (|A|) = \mbox{tr} \, (|A^*|)
\end{equation}
for arbitrary operators $A, \, B \in \cN_1$ and $C \in
L(l^2(\sZ^+))$. For details see \cite{ReS1}, Section VI.6. Let now
$A, \, B \in \mathfrak{F}(\cP)$. Then
\[
\mbox{tr} \, (|P_n(A+B)-(A+B)P_n|) \le \mbox{tr} \, (|P_nA -
AP_n|) + \mbox{tr} \, (|P_nB - BP_n|)
\]
and
\begin{eqnarray*}
\mbox{tr} \, (|P_n(AB)-(AB)P_n|) & = & \mbox{tr} \, (|(P_nA -
AP_n)B + A(P_nB - BP_n)|) \\
& \le & \|B\| \, \mbox{tr} \, (|P_nA - AP_n|) + \|A\| \, \mbox{tr}
\, (|P_nB - BP_n|)
\end{eqnarray*}
by (\ref{e43}) and (\ref{e44}), which implies that $A+B$ and $AB$
are in $\mathfrak{F}(\cP)$ again. Further, if $A_m \in
\mathfrak{F}(\cP)$ and $A_m \to A$ in the norm of $L(l^2(\sZ^+))$,
then
\begin{eqnarray*}
\mbox{tr} \, (|P_nA - AP_n|) & \le & \mbox{tr} \, (|P_n(A-A_m) -
(A-A_m)P_n |) + \mbox{tr} \, (|P_nA_m - A_mP_n|) \\
& \le & 2 \, \mbox{tr} \, P_n \, \|A - A_m\| + \mbox{tr} \,
(|P_nA_m - A_mP_n|),
\end{eqnarray*}
which gives the closedness of $\mathfrak{F}(\cP)$ in
$L(l^2(\sZ^+))$. The symmetry of $\mathfrak{F}(\cP)$ is a
consequence of (\ref{e45}). \hfill \qed \\[3mm]
Recall from Section \ref{ss21} the definitions of the algebra $\cF^\cP$
and of the strong limit homomorphism $W$. Let $\cS(\mathfrak{F}(\cP))$ 
stand for the smallest closed subalgebra of $\cF^\cP$ which contains 
all finite sections sequences $(P_nAP_n)$ where $A$ is in 
$\mathfrak{F}(\cP)$. The following result is the key to several 
generalizations of the first Szeg\"o limit theorem.
\begin{theo} \label{t46}
Let $\bA := (A_n) \in \cS(\mathfrak{F}(\cP))$. Then
\begin{equation} \label{e47}
\frac{1}{n} \mbox{\rm tr} \, (|A_n - P_n W(\bA) P_n|) \to 0
\end{equation}
as $n \to \infty$.
\end{theo}
{\bf Proof.} By (\ref{e44}), the functionals
\[
L(\mbox{im} \, P_n) \to \sC, \quad A_n \mapsto \frac{1}{n} \mbox{tr}
\, (|A_n|)
\]
are uniformly bounded with respect to $n$ (by the constant 1). Hence, 
it is sufficient to prove (\ref{e47}) for sequences $\bA$ in a dense
subalgebra of $\cS(\mathfrak{F}(\cP))$.

Every sequence in $\cS(\mathfrak{F}(\cP))$ can be approximated as
closely as desired (with respect to the norm in $\cF^\cP$) by
sequences of the form
\[
\bB := \sum_j \prod_i (P_n B_{ij} P_n) \quad \mbox{where} \quad
B_{ij} \in \mathfrak{F}(\cP).
\]
Clearly,
\[
W(\bB) = \sum_i \prod_j  B_{ij}.
\]
Thus, and by (\ref{e43}), it is sufficient to prove (\ref{e47})
for sequences of the form $\bB := \prod_i (P_n B_i P_n)$ where
$B_i \in \mathfrak{F}(\cP)$, i.e., to verify that
\begin{equation} \label{e47a}
\frac{1}{n} \mbox{tr} (|P_n B_1 P_n B_2 P_n \ldots P_n B_k P_n -
P_n B_1 B_2 \ldots B_k P_n|) \to 0
\end{equation}
as $n \to \infty$. We prove (\ref{e47a}) in case $k = 2$ from
which the case of general $k$ follows by induction. Assertion
(\ref{e47a}) for $k = 2$ will follow as soon as we have shown
that
\begin{eqnarray*}
\lefteqn{\mbox{tr} \, (|P_n B_1 P_n B_2 P_n - P_n B_1 B_2 P_n|)} \\
&& \le \; \max \, \{\|B_2\| \, \mbox{tr} \, (|P_n B_1 - B_1 P_n|), \,
\|B_1\| \, \mbox{tr} \, (|P_n B_2 - B_2 P_n|) \}
\end{eqnarray*}
for arbitrary operators $B_1, \, B_2 \in L(l^2(\sZ^+))$. This
estimate is a consequence of
\begin{eqnarray*}
\mbox{tr} \, (|P_n B_1 P_n B_2 P_n - P_n B_1 B_2 P_n|)
& = & \mbox{tr} \, (|P_n B_1 (I - P_n) B_2 P_n|) \\
& \le & \|B_1\| \, \mbox{tr} \, (|(I - P_n) B_2 P_n|)
\end{eqnarray*}
and of
\begin{eqnarray*}
\mbox{tr} \, (|(I - P_n) B_2 P_n|) & = & \mbox{tr} \, (|(I-P_n)
(B_2 P_n - P_n B_2)|) \\
& \le & \|I - P_n\| \, \mbox{tr} \, (|P_n B_2 - B_2 P_n|)
\end{eqnarray*}
where we used (\ref{e44}). \hfill \qed \\[3mm]
From (\ref{e42a}) and (\ref{e47}) we conclude that
\[
\frac{1}{n} | \mbox{\rm tr} \, (A_n - P_n W(\bA) P_n)| \to 0.
\]
Thus, if $(w_{ij})_{i, \, j = 0}^\infty$ refers to the matrix 
representation of $W(\bA)$ with respect to the standard basis 
of $l^2(\sZ^+)$, then (\ref{e47}) implies
\begin{equation} \label{e47x}
\left( \frac{\lambda_1 (A_n) + \ldots + \lambda_n(A_n)}{n} -
\frac{w_{00} + \ldots + w_{n-1, \, n-1}}{n} \right) \to 0
\end{equation}
as $n \to \infty$ for every sequence $\bA := (A_n) \in
\cS(\mathfrak{F}(\cP))$.  
\begin{remark}
It is evident that the notion of a F{\o}lner algebra is not restricted 
to the context considered in this section. Indeed, for every sequence 
$\cP = (P_n)$ of orthogonal projections of finite rank acting on a certain 
Hilbert space and tending strongly to the identity operator, there is an 
associated F{\o}lner algebra. This observation allows one to derive  
distributive versions of the first Szeg\"o limit theorem also in the 
higher dimensional context, by employing exactly the same ideas which 
will be pointed out in the following sections. In this way, the results 
of \cite{Lin1,Sio1} can be both easily obtained and generalized. 
\end{remark}
\subsection{Operators and their diagonals} \label{ss42}
A further utilization of (\ref{e47}) and (\ref{e47x}) requires 
to examine the trace $\mbox{tr} \, (P_n W(\bA) P_n)$ which clearly 
depends on the main diagonal of the operator $W(\bA)$ only. In this 
section we show that the main diagonal of operators in 
$\cA_{L^\infty(\sT), \, AP(\sZ)} (\sZ^+)$ behaves quite well.

Let $A \in L(l^2(\sZ))$ be an operator with matrix representation
$(a_{ij})_{i, \, j \in \sZ}$ with respect to the standard basis of
$l^2(\sZ)$. Since
\[
|a_{ii}| = \|P_{\{i\}} A P_{\{i\}}\| \le \|A\|,
\]
the sequence $(a_{ii})_{i \in \sZ}$ belongs to $l^\infty(\sZ)$.
Hence, it defines a multiplication operator on $l^2(\sZ)$ which we
call the {\em main diagonal of} $A$ and which we denote by $D(A)$.
Similarly, the main diagonal of an operator $B \in L(l^2(\sZ^+))$
is defined. It acts as a multiplication operator on $l^2(\sZ^+)$,
and we denote it also by $D(B)$ (which will not rise confusion if
one takes into account where $A$ and $B$ live). In each case,
$\|D(A)\| \le \|A\|$.
\begin{theo} \label{t35}
If $A \in \cA_{L^\infty(\sT), \, AP(\sZ)}(\sZ)$, then $D(A) \in
AP(\sZ)$.
\end{theo}
Of course, then {\em every} diagonal which is parallel to the main
diagonal is almost periodic, too. \\[3mm]
{\bf Proof.} Since $D : L(l^2(\sZ)) \to l^\infty(\sZ)$ is a
continuous linear mapping, and since $AP(\sZ)$ is a closed
subalgebra of $l^\infty(\sZ)$, it is sufficient to prove the
assertion for the case when $A$ is a finite product of Laurent
operators with generating functions in $L^\infty(\sT)$ and of
operators of multiplication by almost periodic functions. Thus, we
can assume that
\[
A = L(a_1) \, b_1 \, L(a_2) \, b_2 \, \ldots \, L(a_k) \, b_k I
\]
with $a_i \in L^\infty (\sT)$ and $b_i \in AP(\sZ)$. Consider the
diagonal $D(A)$ and let $h : \sN \to \sZ$ be an arbitrary
sequence. We have to show that $(U_{-h(n)} D(A) U_{h(n)})_{n \in
\sN}$ has a norm convergent subsequence. Since
\[
U_{-h(n)} D(A) U_{h(n)} = D(U_{-h(n)} A U_{h(n)})
\]
it is sufficient to show that $(U_{-h(n)} A U_{h(n)})_{n \in \sN}$
has a convergent subsequence. Now one has
\begin{eqnarray*}
\lefteqn{U_{-h(n)} A U_{h(n)}} \\
&& = L(a_1) \, (U_{-h(n)} b_1 U_{h(n)}) \, L(a_2) \, (U_{-h(n)}
b_2 U_{h(n)}) \, \ldots  \, L(a_k) \, (U_{-h(n)} b_k U_{h(n)}).
\end{eqnarray*}
Since $b_1$ is almost periodic, there is a subsequence $h_1$ of
$h$ such that the sequence $(U_{-h_1(n)} b_1 U_{h_1(n)})_{n \in
\sN}$ converges. Analogously, there is a subsequence $h_2$ of
$h_1$ such that the sequence $(U_{-h_2(n)} b_2 U_{h_2(n)})_{n \in
\sN}$ converges. We proceed in this way. After $k$ steps we arrive
at a subsequence $g$ of $h$ for which each of the sequences
$(U_{-g(n)} b_i U_{g(n)})_{n \in \sN}$ and, thus, the sequence
$(U_{-g(n)} A U_{g(n)})_{n \in \sN}$ converges. \hfill \qed
\\[3mm]
Let $c_0(\sZ^+)$ stand for the set of all sequences $a : \sZ^+ \to
\sC$ with $a(n) \to 0$ as $n \to \infty$, and write $AP(\sZ^+)$
for the set of all functions $PaP$ where $a \in AP(\sZ)$,
considered as functions on $\sZ^+$. Evidently, both $c_0(\sZ^+)$
and $AP(\sZ^+)$ are closed subalgebras of $l^\infty (\sZ^+)$.
\begin{theo} \label{t36}
If $A \in \cA_{L^\infty(\sT), \, AP(\sZ)}(\sZ^+)$, then $D(A) \in
AP(\sZ^+) + c_0(\sZ^+)$.
\end{theo}
{\bf Proof}. As is the proof of the previous theorem, it is
sufficient to verify the assertion for operators of the form
\begin{eqnarray*}
A & = & T(a_1) \, b_1 \, T(a_2) \, b_2 \, \ldots \, T(a_k) \, b_k
I \\
& = & P L(a_1) P b_1 P L(a_2) P b_2 \, \ldots \, P L(a_k) P b_k P
\end{eqnarray*}
with $a_i \in L^\infty (\sT)$ and $b_i \in AP(\sZ^+)$. We replace
all inner projections $P$ by $I-Q$ and factor out to get
\begin{equation} \label{e37}
A = PBP + R \quad \mbox{where} \quad B \in \cA_{L^\infty(\sT), \,
AP(\sZ)}(\sZ^+)
\end{equation}
and where $R$ is a finite sum, with each item in this sum being a
product of Laurent operators, multiplication operators,
projections $P$ and {\em at least one} projection $Q$. Evidently,
the projections $P$ and $Q$ have a rich operator spectrum, and
$\sigma_+(Q) = \{0\}$. Since the set $L^\$ (l^2(\sZ))$ forms an
algebra we conclude that the operator $R$ has a rich operator
spectrum, too, and the algebraic properties of limit operators
stated in Proposition 1.2.2 in \cite{RRS4} yield that also
$\sigma_+(R) = \{0\}$.

We claim that the main diagonal $D(R) =: \mbox{diag} \, (r_{nn})$
of $R$ is in $c_0(\sZ^+)$. Suppose it is not. Then there is a $C >
0$ and a strongly monotonically increasing sequence $h : \sN \to
\sN$ such that $|r_{h(n),h(n)}| \ge C$ for all $n \in \sN$. Since
$R \in L^\$ (l^2(\sZ))$ there is a subsequence $g$ of $h$ for
which the limit operator $R_g$ exists. Since $h$ (thus, $g$) tends
to $+ \infty$, one has $R_g \in \sigma_+(R)$, whence $R_g = 0$.
This implies in particular that
\[
r_{g(n),g(n)} = P_1 U_{-g(n)} R U_{g(n)} P_1 \to 0,
\]
a contradiction. Thus, $D(R) \in c_0(\sZ^+)$, and passing to the
main diagonals in (\ref{e37}) yields
\[
D(A) = P D(B) P + D(R) \in AP(\sZ^+) + c_0(\sZ^+)
\]
due to Theorem \ref{t35}. \hfill \qed
\begin{prop} \label{p38}
Each function $a \in AP(\sZ^+) + c_0(\sZ^+)$ has a unique
representation in the form $a = PfP + c$ where $f \in AP(\sZ)$ and
$c \in c_0(\sZ^+)$.
\end{prop}
{\bf Proof.} Let $f_1, \, f_2 \in AP(\sZ)$ and $c_1, \, c_2 \in
c_0(\sZ^+)$ be such that $Pf_1P + c_1 = Pf_2P + c_2$. Then $Pf_1P
- Pf_2P = c_2 - c_1$, i.e., $c_2 - c_1 \in c_0(\sZ^+)$ is the
restriction of an almost periodic function. We claim that this 
implies $c_1 = c_2$ and, consequently,  $Pf_1P = Pf_2P$. The latter 
identity further implies $f_1 = f_2$ by Corollary 3.3 in \cite{RRS5}. 

To get the claim, let $f \in AP(\sZ)$ and $c := PfP \in c_0(\sZ^+)$.
Suppose that $c \neq 0$. Then there are an $n_0 \in \sZ^+$ and a 
positive constant $\delta$ with $|c(n_0)| = |f(n_0)| = \delta$. 
Let $h \to + \infty$ be a distinguished sequence for $f$. Then
\begin{eqnarray*}
\|f - U_{-h(n)} f\|_\infty & \ge &
|(f - U_{-h(n)} f)(n_0)| \\
& = & |f(n_0) - f(n_0 + h(n))| \\
& = & |f(n_0) - c(n_0 + h(n))| \to \delta  \quad \mbox{as} \; 
n \to \infty,
\end{eqnarray*}
which is in contradiction to the definition of a distinguished 
sequence.
\hfill \qed
\\[3mm]
Thus, for each operator $A \in \cA_{L^\infty(\sT), \,
AP(\sZ)}(\sZ^+)$, there is a uniquely determined function $f \in
AP(\sZ)$ such that $D(A) - PfP \in c_0(\sZ^+)$. We call this
function the {\em almost periodic part of the main diagonal of}
$A$ and denote it by $D_{ap} (A)$. Note that $D_{ap}(PAP) = D(A)$
for each operator $A \in \cA_{L^\infty(\sT), \, AP(\sZ)}(\sZ)$.
\subsection{The first Szeg\"o limit theorem} \label{ss43}
We are now going to formulate a general version of the first 
Szeg\"o limit theorem which will imply all other versions of Szeg\"o
limit theorems as particular instances. This version is based on a 
fundamental property of every almost periodic function $a$, namely 
that the arithmetic means
\begin{equation} \label{e23}
\frac{1}{n} \sum_{r=0}^{n-1} a(r)
\end{equation}
tend to some value $M(a)$ called the {\em mean value of} $a$ (see
\cite{Cor1}, Theorem 1.28 or \cite{HeR1}, Example (b) in Section
(18.15)). 
\begin{theo} \label{t38a}
Let $\bA = (A_n) \in \cS_{L^\infty(\sT), \, AP(\sZ)}(\sZ^+)$.
Then
\begin{equation} \label{e38b}
\lim_{n \to \infty} \frac{\lambda_1 (A_n) + \dots +
\lambda_n (A_n)}{n} = M(D_{ap}(W(\bA))).
\end{equation}
\end{theo}
{\bf Proof.} It is shown in Corollary 1 in \cite{SeL1} and in
Section 7.2.1 of \cite{HRS2} that the F{\o}lner algebra
$\mathfrak{F}(\cP)$ contains all Laurent operators and all
band-dominated operators. Hence, $\cA_{L^\infty(\sT), \,
AP(\sZ)}(\sZ^+)$ is a subalgebra of the F{\o}lner algebra, and
(\ref{e47}) and (\ref{e42a}) imply
\begin{equation} \label{e38c}
\frac{1}{n} |\mbox{tr} \, (A_n - P_n W(\bA) P_n)|
= \frac{1}{n} |\mbox{tr} \, (A_n) - \mbox{tr} \, (P_n W(\bA) P_n)|
\to 0.
\end{equation}
Evidently, $\mbox{tr} \, (A_n) = \lambda_1 (A_n) + \dots +
\lambda_n (A_n)$, and it remains to show that
\begin{equation} \label{e48}
\frac{1}{n} \mbox{tr} \, (P_n W(\bA) P_n) \to M(D_{ap}(W(\bA))).
\end{equation}
Since $\bA \in \cS_{L^\infty(\sT), \, AP(\sZ)}(\sZ^+)$,
one has $W(\bA) \in \cA_{L^\infty(\sT), \, AP(\sZ)}(\sZ^+)$.
Then, by Proposition \ref{p38},
\[
\frac{1}{n} \mbox{tr} \, (P_n W(\bA) P_n) = \frac{1}{n} \mbox{tr} \,
(P_n D(W(\bA)) P_n) = \frac{1}{n} \left( \sum_{k=1}^n D_{ap}(W(\bA))
(k) + \sum_{k=1}^n c(k) \right)
\]
with a certain function $c \in c_0(\sZ^+)$. Since $\frac{1}{n}
\sum_{k=1}^n c(k) \to 0$, and by what has been said before Theorem
\ref{t38a}, the convergence (\ref{e48}) follows. \hfill \qed \\[3mm]
Note that it is exactly the mean value property of the almost periodic
functions which allows us to prove the existence of the limit in
(\ref{e38b}).
\begin{remark} \label{r48a}
For Toeplitz operators, the block case is considered as being of
particular interest. In order to see how the block case follows from
Theorem \ref{t38a} we mention an obvious generalization of that
theorem. Let $\eta : \sN \to \sN$ be a strongly monotonically
increasing sequence. In place of the sequence $\bA = (A_n) \in
\cS_{L^\infty(\sT), \, AP(\sZ)}(\sZ^+)$ we consider its subsequence
$(A_{\eta(n)})$. Then the limit
\[
\lim_{n \to \infty} \frac{\mbox{tr} \, (A_{\eta(n)})}
{\mbox{tr} \, (P_{\eta(n)})}
= \lim_{n \to \infty} \frac{\lambda_1 (A_{\eta(n)}) + \dots +
\lambda_{\eta(n)} (A_{\eta(n)})}{\eta(n)}
\]
exists and is equal to $M(D_{ap}(W(\bA)))$. The block case follows if
one allows for $d$-periodic coefficients only and if one chooses
$\eta(n) := dn$.
\end{remark}
\section{Special cases} \label{s5}
\subsection{Szeg\"o-type theorems} \label{ss51}
\paragraph{Continuous functions of sequences.}
Here we are going to derive versions of Theorem \ref{t38a} which
hold for functions of sequences in $\cS_{L^\infty(\sT), \,
AP(\sZ)} (\sZ^+)$. Of course, they cannot yield anything which is
substantially new since continuous functions of normal elements of
this algebra belong to $\cS_{L^\infty(\sT), \, AP(\sZ)} (\sZ^+)$
again. But they will bring us closer to the formulation of the
classical Szeg\"o limit theorems.
\begin{theo} \label{t39a}
Let $\bA = (A_n)$ be a normal sequence in $\cS_{L^\infty(\sT), \,
AP(\sZ)}(\sZ^+)$, and let $g$ be any function which is continuous
on a neighborhood in $\sR$ of the stability spectrum $\sigma(\bA + \cG)$.
Then
\begin{equation} \label{e39b}
\lim_{n \to \infty} \frac{g(\lambda_1 (A_n)) + \dots +
g(\lambda_n(A_n))}{n} = M(D_{ap}(g(W(\bA)))).
\end{equation}
\end{theo}
{\bf Proof.} Let $U$ be a neighborhood of $\sigma(\bA + \cG)$ in
$\sR$ and $g$ continuous on $U$. By Proposition \ref{p9b},
\[
\sigma(A_n) \subseteq U \quad \mbox{and} \quad \sigma(W(\bA))
\subseteq U,
\]
and $A_n$ and $W(\bA)$ are normal. Thus, $g(A_n)$ and $g(W(\bA))$
are well-defined via the continuous functional calculus for normal
elements of a $C^*$-algebra (Theorem 6.2.7 in \cite{Aup1}).
Without loss we can also assume that $\sigma_{\cF^\cP} (\bA)
\subseteq U$ such that $g(\bA)$ is well-defined. Indeed, the
spectrum of $\bA$ in $\cF^\cP$ is the union of all spectra
$\sigma(A_n)$ with the stability spectrum of $\bA$. Thus, there is
a finitely supported sequence $\bG$ such that the spectrum of
$(B_n) = \bB := \bA + \bG$ lies in $U$. Since $B_n = A_n$ for
sufficiently large $n$ and since $W(\bB) = W(\bA)$, one can
replace $\bA$ by $\bB$ without loss. Clearly, one also has $B_n =
g(A_n)$ for sufficiently large $n$.

Applying (\ref{e38b}) to the sequence $g(A)$ yields
\begin{equation} \label{e39c}
\lim_{n \to \infty} \frac{\lambda_1 (g(A_n)) + \dots +
\lambda_n (g(A_n))}{n} = M(D_{ap}(W(g(\bA)))).
\end{equation}
The continuous functional calculus for normal elements (or the
Gelfand-Naimark theory for commutative $C^*$-algebras) further
tells us that
\begin{equation} \label{e39d}
\sigma(g(A_n)) = g(\sigma(A_n))
\end{equation}
for all $n$ with $\sigma(A_n) \subseteq U$. Thus,
\begin{equation} \label{e39e}
\lambda_1 (g(A_n)) + \dots + \lambda_n (g(A_n)) =
g(\lambda_1 (A_n)) + \dots + g(\lambda_n (A_n)).
\end{equation}
Finally one has
\begin{equation} \label{e39f}
W(g(\bA)) = g(W(\bA)).
\end{equation}
This equality is evident when $g(\lambda) = p(\lambda, \,
\overline{\lambda})$ where $p$ is a polynomial in two variables,
in which case one has
\[
g(W(\bA)) = p(W(\bA), \, W(\bA)^*),
\]
and it follows for general $g$ since every compactly supported
continuous function can be uniformly approximated by polynomials
of the form $\lambda \mapsto p(\lambda, \, \overline{\lambda})$
due to the Stone-Weierstra{\ss} theorem (Theorem IV.10 in
\cite{ReS1}). The equalities (\ref{e39c}), (\ref{e39e}) and
(\ref{e39f}) imply the assertion. \hfill \qed 
\paragraph{Holomorphic functions of sequences.} 
Next we will discuss a version for non-normal elements which has 
to be based on the holomorphic functional calculus. Recall that, 
for each element $b$ of a Banach algebra $\cB$ with identity $e$ 
and for each function $g$ which is holomorphic in a neighborhood 
$U$ of $\sigma_\cB(b)$, the element $g(b)$ is defined by
\begin{equation} \label{e39g}
g(b) := \frac{1}{2 \pi i} \int_\Gamma g(\zeta) (\zeta e - b)^{-1}
\, d\zeta
\end{equation}
where $\Gamma$ is a smooth oriented Jordan curve in $U \setminus
\sigma_\cB(b)$ which surrounds $\sigma_\cB(b)$. This definition is
independent of the choice of $\Gamma$, and it settles a
homomorphism from the algebra of the holomorphic functions on $U$
into $\cB$ which is continuous in the sense that if a sequence
$(g_n)$ converges to $g$ uniformly on compact subsets of $U$, then
$g(b) = \lim g_n(b)$ in the norm of $\cB$. Moreover,
\begin{equation} \label{e39h}
\sigma_\cB (g(b)) = g(\sigma_\cB (b)).
\end{equation}
For details see \cite{Aup1}, Section III.3.
\begin{theo} \label{t39i}
Let $\bA = (A_n)$ be a sequence in $\cS_{L^\infty(\sT), \,
AP(\sZ)} (\sZ^+)$, and let $g$ be any function which is
holomorphic on a neighborhood $U$ in $\sC$ of the stability
spectrum $\sigma(\bA + \cG)$. Then
\begin{equation} \label{e39j}
\lim_{n \to \infty} \frac{g(\lambda_1 (A_n)) + \dots +
g(\lambda_n(A_n))}{n} = M(D_{ap}(g(W(\bA)))).
\end{equation}
\end{theo}
{\bf Proof.} The proof runs completely parallel to that of Theorem
\ref{t39a}. As there one checks that all occurring terms as well
as the sequence $g(\bA)$ are well defined (the latter after
modification by a finitely supported sequence if necessary). Thus,
the analogue of (\ref{e39c}) holds.

Further, the equality (\ref{e39h}) implies the analogue of
(\ref{e39d}) which, on its hand, yields the analogue of (\ref{e39e}).
Finally, the analogue of (\ref{e39f}) follows by applying the
(continuous and unital) homomorphism $W$ to the contour integral
(\ref{e39g}): approximate this integral by a sequence of Riemann
sums $r_n(\bA)$ and use that $W(r_n(\bA)) = r_n(W(\bA))$. \hfill
\qed \\[3mm]
Another approach to this theorem employs Runge's approximation
theorem (\cite{Gai1}, Theorem 2 in Section III.1) in place of the
holomorphic functional calculus. Runge's theorem yields
approximations of $g(b)$ by linear combinations of $(\zeta_i e -
b)^{-1}$ with simple poles $\zeta_i$ in $U \setminus \sigma(b)$.
(Note that the Riemann sums for (\ref{e39g}) also yield such
approximations.)
\paragraph{Finite sections sequences.}
Next we specify these results to finite sections sequences
$(P_nAP_n)$ where $A$ is a normal operator in $\cA_{L^\infty(\sT),
\, AP(\sZ)}(\sZ^+)$.
\begin{theo} \label{t39}
Let $A$ be a normal operator in $\cA_{L^\infty(\sT), \,
AP(\sZ)}(\sZ^+)$ and let $g$ be any continuous function on the
convex hull of the spectrum of $A$. Then
\begin{equation} \label{e40}
\lim_{n \to \infty} \frac{g(\lambda_1(P_nAP_n)) + \dots +
g(\lambda_n(P_nAP_n))}{n} = M(D_{ap}(g(A))).
\end{equation}
\end{theo}
{\bf Proof.} The interesting new point is that $g$ is merely
assumed to be continuous on the convex hull $I$ of the spectrum of
the operator $A$. Of course, the operator $g(A)$ is still well
defined. Further one knows that all eigenvalues of $P_nAP_n$
belong to $I$, too. This can be most easily seen by introducing
the numerical range
\[
N(B) := \{ \langle Bx, \, x \rangle : x \in l^2(\sZ^+), \, \|x\| = 1 \}
\]
of an operator $B \in L(l^2(\sZ^+))$. It is well known that
\[
\mbox{conv} \, \sigma(A) \subseteq \mbox{clos} \, N(A)
\]
for each operator $A \in L(l^2(\sZ^+))$ and that equality holds in
this inclusion if $A$ is normal (see \cite{BoD1} or Section 3.4.1
in \cite{HRS2}). Here, $\mbox{conv} \, M$ stands for the convex
hull of the set $M \subset \sC$. Consequently, for each normal
operator $A$,
\[
\sigma (P_nAP_n) \subseteq \mbox{clos} \, N(P_n A P_n) \subseteq
\mbox{clos} \, N(A) = \mbox{conv} \, \sigma (A)
\]
where the second inclusion holds since each unit vector in
$\mbox{im} \, P_n$ is also a unit vector in $l^2(\sZ^+)$. Thus,
$g(P_nAP_n)$ is also well-defined. The inclusions $\sigma(P_nAP_n)
\subseteq I$ holding for every $n \in \sN$ together with the
property of being normal further imply that the stability spectrum
of the finite sections sequence $(P_nAP_n)$ is in $I$, too. \hfill
\qed \\[3mm]
In a similar way, one derives the following special case of 
Theorem \ref{t39i}.
\begin{theo} \label{t39k}
Let $A \in \cA_{L^\infty(\sT), \, AP(\sZ)}(\sZ^+)$ and $\bA = 
(P_n A P_n)$. Further, let $g$ be any function which is
holomorphic on a neighborhood $U$ in $\sC$ of the stability
spectrum $\sigma(\bA + \cG)$. Then
\begin{equation} \label{e39l}
\lim_{n \to \infty} \frac{g(\lambda_1(P_nAP_n)) + \dots +  
g(\lambda_n(P_nAP_n))}{n} = M(D_{ap}(g(A))).
\end{equation}
\end{theo}  
Let now $A \in \cA_{AP(\sZ)}(\sZ^+)$ be a band-dominated operator
with almost periodic coefficients. Then we can determine the stability
spectrum of the finite sections sequence $(P_n A P_n)$ by means of
Theorem \ref{t72.5}. If we pass from $(P_n A P_n)$ to a subsequence
$(P_{h(n)} A P_{h(n)})$ then the stability spectrum will decrease in 
accordance with Theorem \ref{t72.6} and, thus, the set of the 
holomorphic functions $g$ for which (\ref{e39l}) holds will become 
larger. The minimal possible stability spectrum (thus, the maximal 
set of holomorphic functions $g$ for which (\ref{e39l}) holds) is 
obtained if we choose $h$ as a distinguished sequence of $A$. In 
this case, the stability spectrum of the sequence 
$(P_{h(n)} A P_{h(n)})$ is equal to 
\[
\sigma(PAP) \cup \sigma(JQAQJ)
\]
by Theorem \ref{t18}.   
\paragraph{Operators in the Toeplitz algebra.} 
Let now $A$ be a normal operator in the Toeplitz algebra 
$\cA_{L^\infty(\sT), \, \sC} (\sZ^+)$ and let $g$ be continuous. 
Then $D_{ap}(g(A))$ coincides with the 0th Fourier coefficient
$g(s_A)_0$ of the function $g(s_A)$ where the symbol $s_A$ of $A$ 
is defined after Proposition \ref{p31}. This equality follows by 
a similar reasoning as in the proofs of Theorems \ref{t35} and 
\ref{t36}. Since $D_{ap}(g(A))$ is a constant function, one 
clearly has $M(D_{ap}(g(A))) = g(s_A)_0$. Thus, specifying Theorem 
\ref{t39} to operators in the Toeplitz algebra yields the following 
version of Szeg\"o's first limit theorem which is due to SeLegue 
\cite{SeL1}.
\begin{coro}[SeLegue] \label{c49}
Let $A$ be a normal operator in $\cA_{L^\infty(\sT), \,
\sC}(\sZ^+)$ and let $g$ be any continuous function on the convex
hull of the spectrum of $A$. Then
\begin{equation} \label{e50}
\lim_{n \to \infty} \frac{g(\lambda_1 (P_nAP_n)) + \dots +
g(\lambda_n (P_nAP_n))}{n} = g(s_A)_0 = \frac{1}{2 \pi} \int_0^{2
\pi} g(s_A(e^{it})) \, dt.
\end{equation}
\end{coro}
In particular, if $A = T(a)$ is a Toeplitz operator with a
generating function $a \in L^\infty(\sT)$, then $s_A = a$. Thus, a
further specification of Corollary \ref{c49} to the case of normal
Toeplitz operators yields the following.
\begin{coro} \label{c49a}
Let $a \in L^\infty(\sT)$ be such that the Toeplitz operator
$T(a)$ is normal, and let $g$ be any continuous function on the
convex hull of the essential range of $a$. Then
\begin{equation} \label{e49b}
\lim_{n \to \infty} \frac{g(\lambda_1 (T_n(a)) + \dots +
g(\lambda_n (T_n(a))}{n} \, = \, \frac{1}{2 \pi}
\int_0^{2 \pi} g(a(e^{it})) \, dt.
\end{equation}
\end{coro}
In this form, one finds the first Szeg\"o theorem in \cite{BSi2},
Theorem 5.10, for instance. Note that a Toeplitz operator $T(a)$
is normal if and only if it is a complex linear combination of a
self-adjoint Toeplitz operator and the identity and, thus, if and
only if the essential range of $a$ is contained in a line segment
(the Brown-Halmos theorem, see Section 3.3 in \cite{BSi2}. Thus,
for Toeplitz operators, there is no basic difference between the
normal and the self-adjoint case. Note also that the finite
sections $P_n T(a) P_n$ are normal for a normal Toeplitz operator.

A final specification of Corollary \ref{c49a} to self-adjoint
Toeplitz operators yields precisely Theorem \ref{t4}. Its
holomorphic version Theorem \ref{t6} follows by a similar
specification of Theorem \ref{t39i}.
\paragraph{Operators in algebras with unique tracial state.}
We finish this section with a few remarks on subalgebras $\cB$ of
the F{\o}lner algebra which own a unique tracial state, i.e., a
state $\tau$ with $\tau(AB) = \tau(BA)$ for each pair of operators
$A, \, B \in \cB$. Their importance for generalized Szeg\"o
theorems rests on the following result. For its proof and all
further facts cited here see \cite{Arv3,Bed1} or Sections 7.2.1
and 7.2.4 in \cite{HRS2}.
\begin{theo}[Arveson, B\'{e}dos]  \label{t51}
Let $\cB$ be a unital $C^*$-subalgebra of the F{\o}l\-ner algebra
${\mathfrak F}(\cP)$. For every $n \ge 1$, let $\rho_n$ be the
state of $\cB$ defined by
\[
\rho_n(A) := \frac{1}{n} \mbox{\rm tr} \, (P_n A P_n),
\]
and let $\cR_n$ be the $^*$-weak-closed convex hull of the set $\{
\rho_n, \, \rho_{n+1}, \, \rho_{n+2}, \dots \}$. Then $\cR_\infty
:= \cap_{n \ge 1} \cR_n$ is a non-empty set of tracial states of
$\cB$.
\end{theo}
Thus, if $\cB$ has a {\em unique} tracial state $\tau$ then the
$\rho_n$ converge $^*$-weakly to $\tau$. In particular,
\[
\lim_{n \to \infty} \rho_n(g(A)) = \tau(g(A))
\]
for each self-adjoint operator $A \in \cB$ and each continuous
function $g$. This implies easily the following version of the
first Szeg\"o limit theorem.
\begin{theo}[Arveson, B\'{e}dos] \label{t52}
Let $\cB$ be a unital $C^*$-subalgebra of the F{\o}l\-ner algebra
${\mathfrak F}(\cP)$ which possesses a unique tracial state
$\tau$. Let further $A \in \cB$ be a self-adjoint operator. Then,
for every compactly supported continuous function $g: \sR \to
\sR$,
\[
\lim_{n \to \infty} \frac{g(\lambda_1(P_nAP_n)) + \dots +
g(\lambda_n(P_nAP_n))}{n} = \tau(g(A)).
\]
\end{theo}
Note that, for each self-adjoint operator $A \in \cB$, the state
$\tau$ gives rise to a natural probability measure $\mu_A$ on
$\sR$ via
\begin{equation} \label{e53}
\int_{-\infty}^\infty g(x) \, d \mu_A(x) := \tau(g(A)).
\end{equation}
A particular example of a $C^*$-subalgebra of the F{\o}lner
algebra with a unique tracial state is the irrational rotation
algebra. The operators in this algebra can be also considered as
band-dominated operators with almost periodic coefficients. Thus,
they are subject both to the Arveson-B\'{e}dos Theorem \ref{t52}
and to our Theorem \ref{t39}. This observation allows one to
identify the tracial state $\tau$ of the irrational rotation
algebra as well as the measures associated with $\tau$ by
(\ref{e53}) via
\[
\int_{-\infty}^\infty g(x) \, d \mu_A(x) = \tau(g(A)) =
M(D_{ap}(g(A))),
\]
which holds for each compactly supported continuous function $g$.
\subsection{Avram-Parter-type theorems} \label{ss52}
The Avram-Parter theorem establishes a formula for the trace of
\[
g(P_n T(\overline{a}) P_n T(a) P_n) \quad \mbox{with} \; a \in
L^\infty(\sT)
\]
and is, thus, immediately related with products of finite sections
sequences and with algebras generated by them. Indeed, we will see 
that this theorem can be considered as another simple special
case of Theorem \ref{t38a}. For each $n \times n$-matrix $B$, let
$\sigma_i(B)$ with $i = 1, \, \ldots, \, n$ refer to the singular
values of $B$, i.e., to the non-negative square roots of the
eigenvalues of $B^*B$. The order of enumeration is again not of
importance.

Let $\bA = (A_n) \in \cF^\cP$. Then the entries of the sequence
$\bB := (\bA^* \bA)^{1/2}$ are the matrices $B_n :=
(A_n^*A_n)^{1/2}$, and
\[
\sigma_j (A_n) = \lambda_j (B_n) \quad \mbox{for} \; j = 1, \,
\ldots, \, n
\]
under suitable enumeration. Thus, application of Theorem
\ref{t39a} to the sequence $\bB$ yields the following.
\begin{theo} \label{t54}
Let $\bA = (A_n)$ be a sequence in $\cS_{L^\infty(\sT), \,
AP(\sZ)}(\sZ^+)$, and let $g$ be any function which is continuous
on a neighborhood in $\sR$ of the stability spectrum $\sigma(\bB +
\cG)$ with $\bB := (\bA^*\bA)^{1/2}$. Then
\begin{equation} \label{e55}
\lim_{n \to \infty} \frac{g(\sigma_1 (A_n)) + \dots +
g(\sigma_n (A_n))}{n} = M(D_{ap}(g(W(\bB)))).
\end{equation}
\end{theo}
\begin{coro} \label{c56}
Let $\bA := (P_nAP_n)$ with $A \in \cA_{L^\infty(\sT), \,
AP(\sZ)}(\sZ^+)$, and let $g$ be any function which is continuous
on a neighborhood in $\sR$ of the stability spectrum $\sigma(\bB +
\cG)$ with $\bB := (\bA^*\bA)^{1/2}$. Then
\begin{equation} \label{e57}
\lim_{n \to \infty} \frac{g(\sigma_1 (A_n)) + \dots +
g(\sigma_n (A_n))}{n} = M(D_{ap}(g(B)))
\end{equation}
with $B := (A^*A)^{1/2}$.
\end{coro}
Further specification to the case of operators in the Toeplitz
algebra yields the following version of SeLegue's result 
(Corollary \ref{c49}).
\begin{coro} \label{c58}
Let $\bA := (P_nAP_n)$ with $A \in \cA_{L^\infty(\sT), \,
\sC}(\sZ^+)$, and let $g$ be any function which is continuous on a
neighborhood in $\sR$ of the stability spectrum $\sigma(\bB +
\cG)$ with $\bB := (\bA^*\bA)^{1/2}$. Then
\begin{equation} \label{e59}
\lim_{n \to \infty} \frac{g(\sigma_1 (A_n)) + \dots +
g(\sigma_n (A_n))}{n} = g(s_B)_0 = \frac{1}{2 \pi} \int_0^{2 \pi}
g(s_B(e^{it})) \, dt
\end{equation}
with $B := (A^*A)^{1/2}$.
\end{coro}
Finally, if $A = T(a)$ is a Toeplitz operator with generating
function $a \in L^\infty(\sT)$, then
\[
s_B = s_{(A^*A)^{1/2}} = (\overline{a} a)^{1/2} = |a|.
\]
\begin{coro}[Avram/Parter] \label{c60}
Let $\bA := (P_n T(a) P_n)$ with $a \in L^\infty(\sT)$, and let
$g$ be any function which is continuous on a neighborhood in $\sR$
of the stability spectrum $\sigma(\bB + \cG)$ with $\bB :=
(\bA^*\bA)^{1/2}$. Then
\begin{equation} \label{e61}
\lim_{n \to \infty} \frac{g(\sigma_1 (A_n)) + \dots +
g(\sigma_n (A_n))}{n} = \frac{1}{2 \pi} \int_0^{2 \pi}
g(|a(e^{it})|) \, dt.
\end{equation}
\end{coro}
This result was established by Parter \cite{Par1} for locally
self-adjoint (= products of continuous and real-valued) generating
functions $a$, and Avram \cite{Avr1} proved it for general
$L^\infty (\sT)$-functions. The algebraic approach to the
Avram/Parter theorem goes back to B\"ottcher and one of the
authors (Section 5.6 in \cite{BSi2}). There (Section 4.5) one also
finds a short illustrated history of the Avram/Parter theorems
which were aimed to explain Moler's phenomenon concerning the
singular value distribution of Toeplitz matrices.

We would also like to mention that Tyrtyshnikov \cite{Tyr1,Tyr2}
was able to show that Corollary \ref{c60} remains valid for
arbitrary functions $a \in L^2(\sT)$ (in which case the Toeplitz
operator $T(a)$ is no longer bounded and our techniques do not
seem to apply).
\subsection{B\"ottcher-Otte-type theorems} \label{ss53}
The continuous and holomorphic functional calculus can also be 
applied to the sequences considered in Theorem \ref{t46} and 
in (\ref{e47x}). It seems that B\"ottcher and Otte \cite{BOt1} 
were interested in results of that type for the first time. The 
following two corollaries to Theorem \ref{t46} follow by a 
straightforward application of the functional calculus as in 
Subsection \ref{ss51}.
\begin{coro} \label{c61a}
Let $\bA = (A_n)$ be a normal sequence in $\cS(\mathfrak{F}(\cP))$, 
and let $g$ be any function which is continuous on a neighborhood 
in $\sR$ of the stability spectrum $\sigma(\bA + \cG)$.
Then
\begin{equation} \label{e61b}
\lim_{n \to \infty} \frac{1}{n} \left( \mbox{\rm tr} \, g(A_n)
- \mbox{\rm tr} \, (P_n g(W(\bA)) P_n) \right) = 0.
\end{equation}
\end{coro}
\begin{coro} \label{c61c}
Let $\bA = (A_n)$ be a sequence in $\cS(\mathfrak{F}(\cP))$, and let 
$g$ be any function which is holomorphic on a neighborhood $U$ in 
$\sC$ of the stability spectrum $\sigma(\bA + \cG)$. Then
$(\ref{e61b})$ holds. 
\end{coro}  
{\small Authors' addresses: \\[3mm]
Steffen Roch, Technische Universit\"at Darmstadt, Fachbereich
Mathematik, Schlossgartenstrasse 7, 64289 Darmstadt,
Germany. \\
E-mail: roch@mathematik.tu-darmstadt.de \\[2mm]
Bernd Silbermann, Technische Universit\"at Chemnitz, Fakult\"at
f\"ur Mathematik, 09107 Chemnitz, Germany. \\
E-mail: bernd.silbermann@mathematik.tu-chemnitz.de}
\end{document}